\documentclass[a4paper,12pt]{article}
\usepackage{amssymb}
\usepackage{amsthm}
\usepackage{amsmath}
\usepackage{latexsym}
\usepackage[dvips]{graphicx}
\usepackage{exscale}
\usepackage[latin1]{inputenc}

\newtheorem{theorem}{Theorem}[section]
\newtheorem{corollary}{Corollary}[section]
\newtheorem{lemma}{Lemma}[section]
\newtheorem{proposition}{Proposition}[section]
\newtheorem{remark}{Remark}[section]
\newtheorem{definition}{Definition}[section]

\newtheorem{assumption}{Assumption}

\newcommand{\ep}{\varepsilon}

\def\ds{\displaystyle}
\def\pa{{\partial}}
\def\na{{\nabla}}
\def\II{\mathcal I}
\def\A{{\mathcal{A}}}
\def\calB{{\mathcal{B}}}
\def\be{\begin{equation}}
\def\ee{\end{equation}}
\newcommand{\RR}{\mathbb R}

\newcommand{\DD}{\mathbb D}
\newcommand{\NN}{\mathbb N}

\newcommand{\RC}{\mathcal R}
\newcommand{\HC}{\mathcal H}
\newcommand{\MM}{\mathcal M}

\newcommand{\beq}{\begin{equation}}
\newcommand{\eeq}{\end{equation}}
\def\beq*{\begin{equation*}}
\def\eeq*{\end{equation*}}
\def\T{{\mathcal{T}}}
\def\E{{\mathbf{E}}}
\def\v{{\mathbf{v}}}
\def\r{{\mathbf{r}}}
\def\X{{\mathbf{X}}}
\def\G{{\mathcal{G}}}
\def\qbar{{\overline{Q}}}
\def\qbarbar{{\overline{\overline{Q}}}}
\def\beqnarray{\begin{eqnarray}}
\def\eeqnarray{\begin{eqnarray}}
\def\beqnarray*{\begin{eqnarray*}}
\def\eeqnarray*{\begin{eqnarray*}}
\def\e{{\mathbf{e}}}
\def\fref#1{(\ref{#1})}
\numberwithin{equation}{section}

\title{Diffusion and guiding center approximation for particle transport in strong magnetic fields
      \thanks{The authors acknowledge support
by the ACI Nouvelles Interfaces des Math\'ematiques No. ACINIM
176-2004 entitled ``MOQUA" and funded by the French ministry of
research, the ACI Jeunes chercheurs no. JC1035 ``Mod\`eles
dispersifs vectoriels pour le transport \`a l'\'echelle
nanom\'etrique" as well as the projet  No.  BLAN07-2 212988 entitled
``QUATRAIN" and funded by the Agence Nationale de la Recherche.}}
\author{Naoufel Ben Abdallah\thanks{Institut de Mathématiques de Toulouse (UMR 5219),
\'equipe MIP\@, Universit\'e Paul Sabatier Toulouse 3,
118 route de Narbonne, 31062 Toulouse Cedex 09, France
(naoufel@\allowbreak math.\allowbreak univ-toulouse.\allowbreak fr)}
\and Raymond El Hajj\thanks{Centre de Mathématiques, Insa de Rennes et IRMAR (UMR 6625),
20 avenue des Buttes de Co\"esmes, 35708 Rennes Cedex 07, France (raymond.el-hajj@insa-rennes.fr).}
}
\date{}

\begin{document}
\maketitle

\begin{abstract}
The diffusion limit of the linear Boltzmann equation with a strong
magnetic field is performed. The giration period of particles around the magnetic field is assumed to be much smaller than the collision relaxation time which is supposed to be much smaller than the macroscopic time. The limiting equation is shown
to be a diffusion equation in the parallel direction while in the
orthogonal direction, the guiding center motion is obtained. The
diffusion constant in the parallel direction is obtained through the
study of a new  collision operator obtained by averaging  the
original one. Moreover, a correction to the guiding center motion is
derived.
\end{abstract}





\section{Introduction}
The motion of charged particles under the action of a strong
magnetic field is an important phenomena encountered in artificial
and natural plasmas (Tokamakas, ionospheric plasmas, etc).
 In the presence of constant
magnetic field $\vec B$, a charged particle (electron in our case)
of mass $m$ and charge $q$ with velocity $v$ follow a helical
trajectory around magnetic field line. The gyration radius of the
particles, called Larmor radius $r_L$, is inversely proportional to
the magnitude of the magnetic field ($\ds r_L=\frac{m v}{q B}$).
This means that when the magnetic field becomes large, the particles
get trapped along the direction of $\vec B$. In addition, when an
electric field $\vec E$ is applied, the particles experience a drift
of the instantaneous centers of Larmor circles, called in the
literature the guiding center motion, in the direction perpendicular
to both the electric and magnetic fields. The velocity of this drift
is given by the following formula
$$v_{drift}=\frac{\vec E \times \vec B}{B^2}.$$
Direct simulations of Vlasov or Boltzmann equations in the presence
of such large magnetic fields requires the numerical resolutions of
the small position and time scales induced by the gyration along the
magnetic field. Hence, the question of deriving approximate models,
numerically less expensive, is of great importance. Various models
have been developed to this aim like the ``guiding-center" and the
gyrokinetic models. The guiding center approximation consists in
averaging the motion over the gyroperiod when supposing $B$ goes to
$+\infty$ ($r_L\rightarrow 0$) \cite{Grad} and the gyrokinetic
approximation is a generalization of the guiding center one in the
case of kinetic equations. For a complete physical review, we refer
for instance to \cite{Lee,Dubin,littlejohn,Northrop}. Mathematically
speaking, homogenization techniques were used in \cite{FS1,FS2} by E.
Fr\'enod and E. Sonnendr\"ucker to justify the high magnetic field
limit of Vlasov and Vlasov-Poisson systems. It is proven in
\cite{FS1} that the tridimensional guiding center approximation
leads to a one-dimensional kinetic model in the direction of the
magnetic field. The drift phenomenon perpendicular to the magnetic
field is rigourously obtained in \cite{FS2} for the two dimensional
Vlasov equation on a sufficiently long time scale. Other asymptotic
regimes of the  Vlasov equation with strong magnetic fields are
studied \cite{FS3,FRS,FW,Frenod,FH}. The gyrokinetic approximation
of the Vlasov-Poisson system has also been intensively studied by F.
Golse and L. Saint-Raymond for the two and tridimensional motions in
various asymptotic regimes
\cite{GLSR1,Saintraymond1,GLSR2,Saintraymond2}.

\medskip

In all the above references, the transport is assumed to be
ballistic, which means that particles do not suffer any collision
during their motion. We are here interested in regimes were
collisions are important and we propose to study the interplay
between the fluid limit induced by collisions and the high frequency
gyrations around the magnetic field. To simplify the study, we shall
only consider collisions with a thermal bath at a given temperature.
In this case, the diffusion limit of the obtained Boltzmann equation
leads to the so-called drift diffusion equation. The question of
deriving diffusions equations from the Boltzmann equation has been
widely studied during the last three decades in the context of
radiative transfer \cite{bensouss,bensouss2,sentis,bar,dautray}, in
semiconductors and plasmas
\cite{poupaud3,GP,MRS,Babd_Deg,Babd_Deg_Gen,Babd_Desv_Gen,degond1,Degond2,Deg_Gou_Pou}
and many other contexts that we do not mention here. In the presence
of a magnetic field, an important quantity to be considered is the
gyroperiod measuring the time that a particle of mass $m$ with
charge $q$ and submitted to the constant magnetic field $B$ makes a
$2\pi$ rotation around this field  ($T_c = 2\pi m/qB$). When the
gyroperiod is much larger than the relaxation time, then the
magnetic field effect disappear in the diffusion limit as will be
seen later on. When $T_c$ is of the same order of magnitude as the
relaxation time, the diffusion matrix has an antisymmetric
component, generated by the magnetic field (this has been proven for
collision with walls by P. Degond \& al \cite{degond1,DLMM,DM}).
Formal results for binary gas mixtures can also be found in
\cite{lucquin,Deg.Lucq.1,Deg.Lucq.2,Degond2}. We shall consider in
this paper the situation where the gyroperiod is much smaller than
the relaxation time. In that case, one expects that  the diffusion
matrix is obtained through the analysis of an averaged collision
operator. Also, the high frequency oscillations only occur in the
direction perpendicular to the magnetic field while parallel
velocity stay unaffected. Therefore, the behaviours of the solutions
in these directions are expected to be different. We shall show in
this paper that the transport is diffusive along the magnetic field,
while it is dominated by the guiding center drift in the orthogonal
direction. The diffusion coefficient in the parallel direction will
be obtained by averaging the collision cross section around the
magnetic field.

\medskip

The outline of the paper is as follows. In the next section, we
introduce the scaling, notations and state the main results. The
proofs rely on the study of the operator accounting for collisions
and gyration around the magnetic field. The analysis of this
operator is done in Section 3. Section 4 is devoted to the
asymptotics of this operator when the scaled gyration period tends
to zero. It is shown that this limit is well described by a
collision operator with a cross section averaged around the magnetic
field. This result is then used in Section 5 in order to prove the
main results of the paper. Some concluding remarks are listed in the
last section of the paper.
\section{Setting of the problem and main results}
The problem we are interested in is a singular perturbation of the
three dimensional Boltzmann equation with a constant large magnetic
field. The position variable is denoted by $\mathbf{r} = (x,y,z)$,
the velocity variable is denoted $\mathbf{v} = (v_x,v_y,v_z)$. The
magnetic field is assumed to be constant and parallel to the $z$
axis. We shall denote by $\mathbf{r}_\perp = (x,y)$ and
$\mathbf{v}_\perp = (v_x,v_y)$, the orthogonal variables. The
distribution function is a solution of the following singularly
perturbed Boltzmann equation
\begin{equation}\label{boltzmann_eps.eta}
\left\{
\begin{array}{ll} &\ds\frac{\partial
f^{\varepsilon\eta}}{\partial t}+ \frac{\T_z f^{\varepsilon\eta}}{\ep}
 +\frac{\T_\perp
f^{\varepsilon\eta}}{\varepsilon\eta}+\frac{(\mathbf{v} \times
 \mathbf{e}_z)}{\varepsilon^2\eta^2} \cdot \nabla_{\mathbf{v}}f^{\varepsilon\eta}=\frac{Q(f^{\varepsilon\eta})}{\varepsilon^2}\\
&f^{\varepsilon\eta}(t=0)=f_0(\mathbf{r},\mathbf{v}).
\end{array}
\right.
\end{equation}
where $\T_\perp$ and $\T_z$ are the respectively perpendicular
and parallel parts of the transport operator
\begin{equation*}
\T_\perp=\v_\perp \cdot \nabla_{\r_\perp} +\E_\perp\cdot \nabla_{\v_\perp},\quad
\T_z=v_z \pa_z
+E_z \pa_{v_z}.
\end{equation*}
The collision operator is given by

\begin{equation}\label{collision}
Q(f)(\v)=\int_{\mathbb{R}^3}\sigma(\v,\v')[\MM(\v)f(\v')-\MM(\v')f(\v)]d\v'.
\end{equation}
The function $\MM$ is the normalized Maxwellian:
\begin{equation}\label{maxwellienne}
\MM(\v)=\frac{1}{(2\pi)^{\frac{3}{2}}}e^{-\frac{1}{2}|\v|^2},
\end{equation}
and the cross section $\sigma$ is symmetric and bounded from below
and above by two positive constants. The parameter $\ep$ is
intended to go to zero and represents the scaled mean free path. The
scaled magnetic field is given by
$\frac{\mathbf{e}_z}{\varepsilon\eta^2}$, while $\mathbf{E} = \left(\begin{array}{c} \mathbf{E}_\perp \\
\mathrm{E}_z\end{array}\right)$ is the scaled electric field. The
gyroperiod is $\ep \eta^2$ and we shall consider the situation
where $\eta$ and $\ep$ go to zero simultaneously and the situation
where $\ep$ tends to zero for $\eta$ given then $\eta$ tends to
zero. Since the magnetic field strongly confines the motion of the
particles in the perpendicular direction $(x,y)$, we have rescaled
this direction with the prefactor $\eta$. This is why the orthogonal
operator ${\mathcal{T}}_\perp$ comes with a different scaling from the
parallel one ${\mathcal{T}}_z$. By doing so, we shall obtain a non
trivial motion in the orthogonal direction (guiding center).

\medskip

\subsection{Scaling}
In this subsection, we shall explain the hypotheses that lead to the scaled Boltzmann equation
(\ref{boltzmann_eps.eta}). In particular, we highlight the differences of the position scales in the parallel and perpendicular directions (with respect to the magnetic field). The starting point is the diffusion scaling in the absence of a magnetic field. The unscaled Boltzmann equation in this case reads
$$\pa_t F + V\cdot \na_X F +{q \na_X \Phi \over m} \cdot \na_V F = \mathbf{Q} (F)$$
where $(t,X,V)$ are the physical time position and velocity variables, while $\Phi$ is the electrostatic potential, $e$ is the elementary charge and $m$ the mass of the electron (to fix the ideas, we consider electrons). The collision operator $\mathbf{Q}$ has the form
$$\mathbf{Q}(F) = \int \Sigma(V,V') [M_T(V) F(V') - M_T(V') F(V)]\, dV',$$
$\Sigma$ is the transition rate and $M_T$ is the Maxwellian at the temperature $T$
$$M_T (V) = ({2\pi K_B T\over m})^{-3/2} \exp (-{mv^2\over 2 K_B T}),$$
$K_B$ being the Boltzmann constant.
The thermal energy and velocity are respectively defined by
$$U_{th} = K_B T, \quad V_{th} = \sqrt{K_B T\over m}.$$
The electrostatic energy is assumed to be of the order of the thermal energy and varies over a macroscopic length scale $L$ :
$$q \Phi(X) = U_{th} \phi({X\over L}).$$
The transition rate $\Sigma$ is the inverse of a time. The typical value of this time, the scattering time, is denoted
by $\tau_s$. So we write
$$\Sigma (V,V')= ={1\over \tau_s} \sigma ({V\over V_{th}}, {V'\over V_{th}})$$
where $\sigma$ is its scaled version assumed to be of order $1$.
The mean free path $L_s$ will be defined as the distance that a thermal electron crosses during the scattering time $\tau_s$
$$L_s = V_{th} \tau_s.$$
This distance is assumed to be much smaller than the macroscopic distance $L$ and we shall denote it by
$$\ep = {L_s\over L}.$$
The only thing left to know is the macroscopic time scale. Since the kernel of the collision operator is generated by the Maxwellian which carries no flux, one has to go to the diffusion scaling and define the macroscopic time $\tau_m$ by
$${\tau_s\over \tau_m} =({L_s\over L})^2 = \ep^2.$$
Using $(\tau_m, L, V_{th}, U_{th})$ as units for time, position, velocity and energy, the Boltzman  equation can be written
$$\pa_t f +{1\over \ep} (\v \cdot \na_\r f - \na_\r \phi \cdot \na_\v f) = {Q(f) \over \ep^2},$$
where $Q$ is the scaled Boltzmann operator (\ref{collision}).

\medskip

Let us now look at the effect of a large magnetic field ${\calB} = B \e_z.$ The scaled Boltzmann equation contains the additional term $-{q \over m} (V\times \calB) \cdot \na_V F.$ An electron submitted to this magnetic field spins around it
following the equation
$$m{d V\over dt} = - q V \times {\calB}.$$
The period of this precession, called the cyclotron period is given by
$${\tau_c }={m \over q B}.$$
We shall assume that this period is much smaller than the scattering time and we write
$$\tau_c = \eta^2 \tau_s.$$
Of course the motion along the magnetic field stays unchanged and we shall not change the length scale in this direction. In the orthogonal direction however, the combined action of the electric field and the magnetic field results in the drift of the
guiding center (the electron motion is a rotation around a center which is now drifting under the action of the electric field). The velocity of this drift is given by
$$V_c = {E \times \calB \over B^2}.$$
Its order of magnitude is then given by
$$V_c ={E_\perp \over B},$$
wher $E_\perp$ is the perpendicular component of the electric field.

\medskip

The last hypothesis that we make is the potential $\Phi$ varies in the orthogonal direction on a lengthscale $L_\perp$
which is not equal to the parallel lengthscale $L$. Moreover, $L_\perp$ is such that the guiding center moves during the macroscopic time by a distance of the order of $L_\perp$ itself. In other words, we want
$$E_\perp = {U_{th} \over q L_\perp}, \quad V_c \tau_m = L_\perp.$$
Replacing the thermal energy by its expression given above in this section, we find after straightforward computations that  $$L_\perp^2 = {\tau_m U_{th} \over q B} = \eta^2L^2.$$
Now rescaling the orthogonal position variable by $L_\perp$, we obtain the Boltzmann equation (\ref{boltzmann_eps.eta}).

\medskip

\begin{remark}{\rm
From the scaling hypotheses, we see that between two successive collisions,  the electron precesses a large number of times around the magnetic field.
Therefore, one might expect that the collision cross section will be averaged along these precessions. On the other hand, since the orthogonal lengthscale has been designed in order to see the guiding center drift, the limiting
equation should exhibit this term. Finally, since the parallel motion to the magnetic field does not feel the magnetic field, one is entitled to expect a non vanishing diffusion in this direction. This is exactly what we shall prove in this paper : drift in the orthogonal direction, finite diffusion in the parallel direction governed by a collision operator averaged over precessions around the magnetic field.}
\end{remark}

\medskip

\subsection{Notations}
 The analysis of the limit $\ep \to 0$ naturally leads to the study of the operator
\be
\label{qeta}
Q^\eta = Q -  \frac{(\mathbf{v} \times
\mathbf{e}_z)}{\eta^2} \cdot \nabla_{\mathbf{v}},
\ee
while  the limit $\eta\to 0$, involves the cylindrical averaging around the vector $\e_z$.
To this aim, we denote by
$\RC(\tau)$  the rotation around $\e_z$ with angle
$\tau$ which is the multiplication operator by the following matrix
\begin{equation}\label{R(tau)}
\RC(\tau)=\left(
          \begin{array}{ccc}

                \cos\tau & \sin\tau & 0 \\
                -\sin\tau & \cos\tau & 0\\
                0 &0 & 1 \\
          \end{array}
        \right).
\end{equation}

 The cylindrical average of a function $h = h(\v)$ is defined by
\be \label{average} \A(h)(\v) := \overline{h}(\v) = \frac{1}{2 \pi}
\int_0^{2\pi} h(\RC(\tau) \v)\, d\tau. \ee It is clear that $\A$ is
the orthogonal projector (in the $L^2$ sense) on the set of
cylindrically invariant functions. We also define the partial
average \be \label{fbartau} \A_\tau(h)(\v) := \overline{h}_\tau (\v)
= {1\over 2 \pi} \int_0^{\tau} h(\RC(s) \v)\, ds. \ee We shall
moreover define the averaged  function $\overline{\sigma}$ and the
symmetric average $\overline{\overline{\sigma}}$ by
$$\overline{\sigma} (\v,\v') = {1\over 2 \pi} \int_0^{2\pi} \sigma\left(\RC(\tau) \v,\v'\right)\, d\tau$$
$$\overline{\overline{\sigma}} \left(\v,\v'\right) = {1\over 4 \pi^2}
\int_0^{2\pi}\int_0^{2\pi} \sigma\left(\RC(\tau) \v,\RC(\tau') \v'\right)\, d\tau d\tau' = {1\over 2 \pi}
\int_0^{2\pi} \overline{\sigma}\left(\v, \RC(\tau')\v'\right)\, d\tau',$$
and the corresponding collision operators
\be
\label{qbar}
\overline{Q}(f)(\v)=\int_{\mathbb{R}^3}\overline{\sigma}(\v,\v')[\MM(\v)f(\v')-\MM(\v')f(\v)]d\v',
\ee
\be
\label{qbarbar}
\overline{\overline{Q}}(f)(\v)=\int_{\mathbb{R}^3}\overline{\overline{\sigma}}(\v,\v')[\MM(\v)f(\v')-\MM(\v')f(\v)]d\v'.
\ee
We define the functional spaces
\begin{equation}\label{L2M}
L^2_\MM=\left\{f=f(\v)\,/\, \int_{\RR^3}\frac{|f(\v)|^2}{\MM(\v)}d\v
 < +\infty \right\},
\end{equation}
and its cylindrically symmetric subspace
\begin{equation}\label{barL^2_MM}
{\overline{L}}^2_\MM=\{f\in L^2_\MM/\quad f(\RC(\tau)\v)=f(\v),\quad  \forall \tau\in[0,2\pi]\}.
\end{equation}
We also define the space
\begin{equation}
\label{h2c}
\HC^2_\MM = \{ f\in L^2_\MM, \quad\mbox{such that},\,\, \pa^\alpha f \in L^2_\MM, \quad \mbox{for } |\alpha|\leq 2\}
\end{equation}
We shall make the following hypotheses
\begin{assumption}\label{crosssectionAssump.}
The cross-section $\sigma$ belongs to $W^{2,\infty}(\RR^6)$ and is
supposed to be symmetric and bounded from above and below:
\begin{equation}\label{crosssection}
\exists \alpha_1,\alpha_2>0,\quad 0<\alpha_1\leq
\sigma(\v,\v')= \sigma(\v',\v)\leq\alpha_2,\quad \forall (\v,\v')\in \RR^6.
\end{equation}
\end{assumption}

\begin{assumption}\label{pot.Assump.}
There exists a potential $V(t,\r)$ in
$C^1(\RR^+;W^{1,\infty}(\RR^3))$, such that $\mathbf{E} = - \nabla_\r V$.
\end{assumption}
We finally define the space
\begin{equation}\label{L2MV}
L^2_{\MM_V}=\left\{f=f(\r,\v)\,/\, \int_{\RR^6}\frac{|f(\r,\v)|^2}{\MM(\v)\mathrm{e}^{-V(\r)}}d\v
 < +\infty \right\}.
\end{equation}

\subsection{Main results}

\noindent The main results of this paper are summarized in the
following two theorems. The first one deals with the limit
$\ep \to 0$ while $\eta > 0$ is kept fixed.

\begin{theorem}\label{main_th.2}
Let $\eta > 0$ be fixed.
Let $T\in\RR^*_+$ and assume that the initial data of
\eqref{boltzmann_eps.eta}, $f_0$, belongs to $L^2_{\MM_V}$.
Then, with Assumptions \ref{crosssectionAssump.} and
\ref{pot.Assump.},
 the problem \eqref{boltzmann_eps.eta} has a unique weak solution in $C^0([0,T], L^2_{\MM_V}).$ Moreover, the
sequence $(f^{\ep\eta})_\ep$
 converges weakly to $\rho_\eta(t,\mathbf{r})\MM(\mathbf{v})$ in
$L^\infty((0,T);L^2_{\MM_V}))$ weak * where  $\rho_\eta$
satisfies
\begin{equation}\label{continuity_eta}
\partial_t\rho_\eta-\mathrm{div}(\DD^\eta(\nabla \rho_\eta-\rho_\eta \E))=0
\end{equation}
and the diffusion matrix $\DD^\eta$ is given by the formula
$$\DD^\eta = \int \left(\begin{array}{c} {1\over \eta} \v_\perp \\ v_z
\end{array}\right)\otimes\left( \begin{array}{c} \eta \X^\eta_\perp \\ X^\eta_z
\end{array}\right)   d\v,$$
$\X^\eta = \left(\begin{array}{c} \X^\eta_\perp \\[3mm] X^\eta_z \end{array}\right)$ being the only solution of (\ref{chi_eta}).
\end{theorem}
\begin{remark}{
\rm In the relaxation time case $\sigma(\v,\v') = {1\over \tau}$, the matrix $\DD^\eta$ can be computed explicitly:
\begin{equation*}
\DD^\eta=\tau \left(\begin{array}{ccc}
                           \frac{\eta^2}{\tau^2+\eta^4} &  \frac{\tau}{\tau^2+\eta^4}& 0 \\
                           -\frac{\tau}{\tau^2+\eta^4} &\frac{\eta^2}{\tau^2+\eta^4} & 0 \\
                           0 & 0 & 1 \\
                                \end{array}
                              \right).
\end{equation*}
In the limit $\eta \to 0$, the upper 2 x 2 block of the matrix
reduces to the antisymmetric matrix ${\mathcal{I}} =
\left(\begin{array}{cc}
0 & 1\\[3mm]
-1& 0
\end{array}\right)$, whereas its  symmetric part is of order
$\eta^2$. The following proposition shows that these features are
still valid in the case of a non constant scattering section
$\sigma$.
}
\end{remark}
\begin{proposition}
\label{chieta} Writing $\mathbf{X}^\eta = \left(\begin{array}{c}\ds
\mathbf{X}^\eta_\perp \\ \mathrm{X}^\eta_z\end{array}\right)$, we
have the following expansions in ${\mathcal{H}}_\MM^2$ as $\eta$ tends to
zero.
$$\mathrm{X}^\eta_z = \mathrm{X}_z^{(0)} + \eta^2 \mathrm{X}_z^{(1)} + {\mathcal{O}}(\eta^4),$$
$$\mathbf{X}^\eta_\perp = \mathbf{X}_\perp^{(0)}+ \eta^2 \mathbf{X}_\perp^{(1)} + {\mathcal{O}}(\eta^4).$$
The diffusion matrix $\DD^\eta$ is a definite positive matrix. Its
symmetric and antisymmetric parts, $\DD^{\eta}_s$ and
$\DD^{\eta}_{as}$, have respectively the following expansions
$$\DD^{\eta}_s = \left(\begin{array}{cc}\ds \eta^2 \int \v_\perp \otimes \X^{(1)}_\perp d\v & 0\\[3mm]
0 & D_z\end{array} \right)  +\eta^3\DD^{(1)}_{z\perp}+{\mathcal{O}}(\eta^4),$$
$$\DD^{\eta}_{as} = \left(\begin{array}{ccc}0&1 & 0\\[3mm]
-1&0 & 0\\[3mm]
0 & 0&0\end{array}\right) + \eta \DD^{(0)}_{z\perp} +{\mathcal{O}}(\eta^4)
$$
where $D_z = \int v_z X^{(0)}_z\, d\v$ and
$$
\DD^{(0)}_{z\perp} = \int\left(\begin{array}{c} \mathbf{0} \\ v_z
\end{array}\right) \otimes  \left( \begin{array}{c} \X^{(0)}_\perp \\ 0
\end{array}\right)  d\v
- \int \left(\begin{array}{c}  \X^{(0)}_\perp \\ 0
\end{array}\right)\otimes\left( \begin{array}{c} \mathbf{0} \\ v_z
\end{array}\right)   d\v
$$
$$
\DD^{(1)}_{z\perp} =  \int \left( \begin{array}{c} \X^{(1)}_\perp \\
0
\end{array}\right) \otimes \left(\begin{array}{c} \mathbf{0} \\ v_z
\end{array}\right) d\v
+ \int \left( \begin{array}{c} \mathbf{0} \\ v_z
\end{array}\right) \otimes \left(\begin{array}{c}  \X^{(1)}_\perp \\ 0
\end{array}\right) d\v.
$$
\end{proposition}
\begin{remark}{\rm
The above theorem  shows that the diffusion is of order ${\mathcal{O}}(1)$ in the parallel direction while it scales like $\eta^2$ in
the orthogonal one, exactly like in the relaxation time coefficient.
The antisymmetric part of the diffusion matrix acts essentially on
the orthogonal direction and leads to the guiding center motion.
However, the above formula shows a correction of order $\eta$ of the
guiding center motion, which is due to collisions. Indeed, it is
readily seen that for any three dimensional vector $\mathbf{Z}$
$$- \DD^{\eta}_{as} \mathbf{Z} = \mathbf{u}^\eta \times \mathbf{Z} +{\mathcal{O}}(\eta^4)$$
where
$$\mathbf{u}^\eta = \left(\begin{array}{c}
- \eta \int X^{(0)}_y v_z d\v \\[3mm] \eta \int X^{(0)}_x v_z d\v \\[3mm] 1
\end{array}\right),$$
and $(X^{(0)}_x, X^{(0)}_y)$ are the components of $\X^{(0)}_\perp$.
This leads to
$$- \nabla \cdot (\DD^{\eta}_{as}
(\nabla_\r \rho_\eta - \rho_\eta \mathbf{E})) =  \nabla \cdot
\left(\rho_\eta E \times \mathbf{u}^\eta\right) + {\mathcal{O}}
(\eta^4).$$
}
\end{remark}

\begin{theorem}\label{main_th.1}
Under the same hypotheses as for Theorem \ref{main_th.2}, we now let
$\ep$ and $\eta$ simultaneously and independently tend to zero (there is no assumption on the relative scale between $\ep$ and $\eta$).  In this case, the sequence
$(f^{\ep\eta})_{\ep\eta}$ of weak solutions of
\eqref{boltzmann_eps.eta} converges weakly to $\rho(t,\r)\MM(\v)$ in
$L^\infty([0,T],L^2_{\MM_V}) $ weak *, where $\rho$ is the solution
of
\begin{equation}\label{centerguide}
\left\{\begin{array}{l}\partial_t \rho + \pa_z \mathrm{J}_z + \nabla \cdot (\rho {\E\times
\e_z} ) =0\\\ds
 \rho(0,\r)=\int_{\RR^3}f_0(\r,\v)\, d\v
\end{array}
\right.
\end{equation}
The  parallel current density, $\mathrm{J}_z$ is given by
\begin{equation*}
\mathrm{J}_z=-D_z(\pa_z \rho - E_z \rho)
\end{equation*}
where $D_z$ is the diffusion constant given by
\begin{equation}\label{difparallel}
D_z=\int_{\RR^3}\mathrm{X}_z^{(0)}  v_z d\v,
\end{equation}
and $\mathrm{X}_z^{(0)}$ is the zero-th order term of $X_z^\eta$  (defined by  \eqref{chi_z^0}).
\end{theorem}


\section{Analysis of the operator $Q^\eta$.}\label{section.r_eta}

\noindent Let us consider the space $L^2_\MM$ introduced in \fref{L2M} and define the scalar product on this space by
\begin{equation*}
\langle f,g\rangle_\MM=\int_{\RR^3}\frac{f  g}{\MM}d\v.
\end{equation*}
Let $D (Q^\eta)$ be defined by
\begin{equation}\label{defReta}
D = D(Q^\eta):=\{f\in L^2_\MM/\,(\v\times\e_z)\cdot \nabla_\v f\in
L^2_\MM\}.
\end{equation}
We have the following result
\begin{proposition}\label{properR}
The operator $Q^\eta$ given by \eqref{qeta} with domain $D$ defined
by \fref{defReta} satisfies the following properties

\begin{enumerate}
\item For any $f\in D$, we have
$\ds \int_{\RR^3} Q^\eta   (f)d\v=0$ (mass conservation).
\item $(-Q^\eta,D)$ is a maximal monotone operator.
\item The kernel of $Q^\eta$ is the real line spanned by the Maxwellian :
\begin{equation*}
Ker(Q^\eta)=\{f\in L^2_\MM/\,\exists \alpha\in\RR\mbox{ such that }
f(\v) =\alpha \MM(\v)\}.
\end{equation*}
\item  Let $\mathcal{P}$ be the orthogonal projection on $Ker
Q^\eta$. The following  coercivity inequality holds for any function $f\in D$
\begin{equation}\label{coercive}
-\langle
Q^\eta(f),f\rangle_\MM\geq\alpha_1\|f-\mathcal{P}(f)\|_\MM^2.
\end{equation}
where $\alpha_1>0$ is the lower bound of $\sigma$
\eqref{crosssection}.
\item The range of $Q^\eta$, denoted by $Im(Q^\eta)$, is the set of functions $g\in L^2_\MM$
satisfying the solvability condition: $\ds \int_{\RR^3}g(\v)d\v=0$.
In addition, for all $g\in Im(Q^\eta)$ there exists a unique
function $f\in D$ satisfying both $Q^\eta(f)=g$ and  $\ds \int_{\RR^3}f(\v)d\v=0$.
\end{enumerate}
\end{proposition}

\begin{proof}
First of all, we recall that the operator $Q$ satisfies all items of
Proposition \ref{properR} (by substituting $D$ by $L^2_\MM$). Now
Items  1., 3. and 4. of the Proposition follow immediately from
these properties since we have $\int_{\RR^3}(\e_z\times \v)\cdot
\nabla_\v(f)d\v=\int_{\RR^3}\mathrm{div_\v}[(\e_z\times \v)f]d\v=0$
and $\langle(\e_z\times v)\cdot \nabla_v f,f\rangle_\MM=0$. In order
to prove Item 2., we need to prove that $I - Q^\eta$ is one-to-one
from $D$ to $L^2_\MM$. In order to do so, we rewrite the equation
$f- Q^\eta (f) = g$ under the following form
\begin{equation*}
\frac{1}{\eta^2}(\v\times\e_z )\cdot \nabla_\v f + (1 +\nu (\v)) f =
Q^+(f) + g,
\end{equation*}
where we have written
\be
\label{collision1}
Q(f) = Q^+(f) - \nu f,
\ee
\begin{equation}\label{K}
Q^+(f)=(\int_{\RR^3}\sigma(\v,\v')f(\v')d\v') \MM(\v)
\end{equation}
\begin{equation}\label{nu}
\nu(\v)=\int_{\RR^3}\sigma(\v,\v')\MM(\v')d\v'.
\end{equation}
The  characteristics of the equation are defined by
{\begin{equation*} \left\{\begin{array}{rl}
\ds \frac{d\v}{dt}&= \ds \frac{1}{\eta^2}(\v\times\e_z)\\
\v|_{t=0}&=\v_0
\end{array}\right.
\end{equation*}}
whose solution is $\ds \v(t)=\RC(\frac{t}{\eta^2})\v_0$. The
solution of the equation can then be defined by integration over the
characteristics and finally leads to $f =  {\mathcal{F}}(Q^+(f) + g)$
where

\begin{align*}{\mathcal{F}} (h) (\v) =&\frac{\eta^2}{
B_\eta(\v)}\int_0^{2\pi}h(\RC(\tau)\v) {e^{\eta^2\int_0^{\tau}(1
+\nu(\RC(s)v))ds}}d\tau, \\  B_\eta(v)=&e^{\eta^2
\int_0^{2\pi}(1+\nu(\RC(s)v))\,ds}-1.
\end{align*}
A fixed point argument leads to the existence and uniqueness of
$f\in L^2_\MM$ solving the equation $f- Q^\eta (f) = g$ in the
distributional sense. Since $Q$ is a bounded operator on $L^2_\MM$,
we immediately obtain that $(\v\times\e_z)\cdot \nabla_\v f \in
L^2_\MM$ so  that the constructed solution is in the domain $D$ of
the operator $Q^\eta$. Let us now prove item 5. Consider $g\in
L^2_\MM$ such that $\ds \int_{\RR^3}g(v)d\v=0$. From the maximal
monotonicity of $Q^\eta$, for each $\lambda>0$, $\lambda Id-Q^\eta$
is surjective. Therefore,  for all $\lambda >0$, there exists
$f_\lambda\in D(Q^\eta)$ satisfying
\begin{equation}\label{f_lambda}
\lambda f_\lambda-Q^\eta (f_\lambda)=g.
\end{equation}

\noindent Integrating this equation with respect to $\v$ gives $\ds
\lambda\int_{\RR^3} f_\lambda d\v=\int_{\RR^3} g d\v=0$ and with
\eqref{coercive} we get $-\langle
Q^\eta(f_\lambda),f_\lambda\rangle_\MM\geq
\alpha_1\|f_\lambda\|_\MM^2$. Taking the scalar product of
\eqref{f_lambda} with $f_\lambda$ in $L^2_\MM$, one obtains

\begin{equation*}
\lambda \|f_\lambda\|_\MM+\alpha_1\|f_\lambda\|_\MM\leq
\|g\|_\MM.
\end{equation*}
\noindent This implies that $(f_\lambda)_{\lambda >0}$ is bounded in
$L^2_\MM$. Up to an extraction of a subsequence, there exists $f\in
L^2_\MM$ such that $(f_\lambda)$ converges weakly to $f$. By passing
to the limit $\lambda \rightarrow 0$ in the weak form of equation
\eqref{f_lambda}, we get
\begin{equation*}
\frac{1}{\eta^2}(\v\times\e_z) \cdot \nabla_v f-Q(f)=g\quad \mbox{in
} D'(\RR^3).
\end{equation*}
In addition, $Q(f)$ and $g$ belong to $L^2_\MM$, then
$(\v\times\e_z)\cdot \nabla_\v f\in L^2_\MM$ and then $f\in
D(Q^\eta)$. Of course $\int_{\RR^3} f d\v = 0$ since this is true
for $f_\lambda$. Now if there is another solution of $-Q(\cdot) = g$
with zero velocity average, the difference of this solution with $f$
is in the kernel of $Q^\eta$ and has a zero average. It is
necessarily equal to zero.

\end{proof}

The following corollary is a direct consequence of Proposition
\ref{properR}.
\begin{corollary}
\label{chichi}
For all $\eta>0$, there exists a unique $\X_\eta = \left(\begin{array}{c} \X^\eta_\perp \\[3mm] \mathrm{X}^\eta_z
\end{array}\right)$ in $(D(Q^\eta))^3$
satisfying
\begin{equation}\label{chi_eta}
-Q^\eta(\X^\eta)=\left(\begin{array}{c} \ds {1\over
\eta^2}\v_\perp\\ v_z\end{array}\right)\MM(\v)\quad \mbox{and }\quad
\int_{\RR^3}\X^\eta(\v)d\v=0
\end{equation}
and from the coercivity inequality \eqref{coercive},
$(\mathrm{X}_z^\eta)_\eta$ and $(\eta^2\X_\perp^\eta)_\eta$ are bounded sequences in $L^2_\MM$. Namely,
we have  the following a priori estimates
\begin{equation}\label{chi_eta_bound} \|\mathrm{X}_z^\eta\|_\MM +  \|\eta^2\X_\perp^\eta\|_\MM\leq
\frac{1}{\alpha_1}\|\v \MM\|_\MM.
\end{equation}
\end{corollary}

The following proposition allows to reformulate the problem $Q^\eta(f) = - g$ by using the characteristics.

\begin{proposition}
For all $g\in Im(Q^\eta)$, we have the following equivalence
\begin{equation}\label{characIm(R_eta)}
-Q^\eta(f)=g \Leftrightarrow f -L_\eta (Q^+(f))=L_\eta(g)
\end{equation}

\noindent where $Q^+$ is given by \eqref{K} and $L_\eta$ is the
operator on $L^2_\MM$ given by the expression:
\begin{equation}\label{L_eta}
L_\eta(f)=\frac{\eta^2}{
C_\eta(v)}\int_0^{2\pi}f(\RC(\tau)v)
{e^{\eta^2\int_0^{\tau}\nu(\RC(s)v)ds}}d\tau, \quad
 {
C_\eta(v)=e^{\eta^2 \int_0^{2\pi}\nu(\RC(s)v)ds}-1. }
\end{equation}
\end{proposition}
\begin{proof}
Let $(S^\eta, D(S^\eta))$ be the unbounded operator on $L^2_\MM$
defined by
\begin{eqnarray}
D(S^\eta)&=&\{f\in L^2_\MM \,/\, (\v\times\e_z )\cdot \nabla_\v f\in
L^2_\MM\} \\ S^\eta(f)&=&\frac{1}{\eta^2}(\v\times\e_z )\cdot
\nabla_\v f+\nu(\v) f \quad \forall f\in D(S^\eta).
\end{eqnarray}

\noindent With the decomposition \eqref{collision1}, we have
$Q^\eta=Q^+-S^\eta$. Proceeding as in the proof of Proposition
\ref{properR} (replacing $\nu + 1$ by $\nu$), we show that $S^\eta$
is invertible and its inverse is nothing but $L_\eta$. The equation
$-Q^\eta(f)=g$ is equivalent to $f
-(S^\eta)^{-1}(Q^+(f))=(S^\eta)^{-1}(g)$ which concludes the proof.
\end{proof}

\begin{lemma}
For all $f\in L^1_{loc}(\RR^3)$, we have:
\begin{equation}\label{maj.L_eta}
|L_\eta(f)|\leq \frac{\eta^2}{e^{2\pi\alpha_1\eta^2}-1}
\int_0^{2\pi}|f(\RC(\tau)v)|e^{\alpha_2\eta^2\tau} d\tau.
\end{equation}

\end{lemma}
\begin{proof}
This estimate follows immediately from the definition of $L_\eta$
\eqref{L_eta} using Assumption \ref{crosssectionAssump.}.
\end{proof}

\begin{proposition}\label{chieta_regular.}
Under Assumption \ref{crosssectionAssump.}, the solution $\X^{\eta}$
of \eqref{chi_eta} belongs to \break$(\HC^2_\MM(\RR^3))^3$ and
$\X^\eta$ and all its derivatives with respect to $v$ of order less
than or equal $2$ (up to dividing by the Maxwellian) are
polynomially increasing when $|v|$ goes to $+\infty$. Namely, we
have
\begin{equation}\label{polyn.increas.chi_eta_1}
|\X^\eta(\v)| \leq Q_1^\eta(|\v|) \MM(\v)
\end{equation}
\begin{equation}\label{polyn.increas.chi_eta_2}
{|\pa^\alpha \X^\eta(\v)|}\leq Q_3^\eta(|\v|){\MM(\v)}
\end{equation}
for all $\alpha=(\alpha_1,\alpha_2)\in \NN^2$ with
$\alpha_1+\alpha_2\leq 2$ and where $\ds
\pa^\alpha:=\sum_{i,j=1}^3\partial_{\v_i}^{\alpha_1}
\partial_{\v_j}^{\alpha_2}$. Here, $Q_1^\eta$ and $Q_3^\eta$ are two polynomials
of degrees $1$ and $3$ respectively (depending on $\eta$).
\end{proposition}

\begin{proof}
 From \eqref{characIm(R_eta)}, we deduce the identity
\begin{equation*}
\X^\eta=L_\eta(\v^\eta\MM)+L_\eta(Q^+(\X^\eta)),
\end{equation*}
where
$$\v^\eta = \left(\begin{array}{c} \ds {1\over \eta^2}\v_\perp\\ v_z\end{array}\right).$$
Moreover, we have $\ds |Q^+(\X^\eta)|\leq
\alpha_2(\int_{\RR^3}\X^\eta(\v')d\v') \MM(\v)\leq c_0 \MM(\v)$ and
with \eqref{maj.L_eta} one obtains \eqref{polyn.increas.chi_eta_1}.
Estimate \eqref{polyn.increas.chi_eta_2} follows by derivation of
the above equation with respect to $v$ and by using the fact that
the cross-section $\sigma(v,v')$ belongs to $W^{2,\infty}(\RR^6)$.

\end{proof}

\section{Expansion of $\X^\eta$ with respect to
$\eta$}\label{expansionchieta}

A rigourous expansion of $\X_\eta$ around $\eta=0$ will be carried
out in this section. Corollary \ref{chichi} provides us with bounds of
order ${\mathcal{O}}(1)$ for $\mathrm{X}_z^\eta$ and of order ${\mathcal{O}}({1\over \eta^2})$ for $\mathrm{X}_\perp^\eta$. We shall prove in
this section that $\X_\perp^\eta$ is actually bounded in
$\HC^2_\MM$. This is due to the fast precession around the axis
$\e_z$ generated by the differential term ${1\over
\eta^2}(\v\times\e_z )\cdot \nabla_\v $. Filtering out these
oscillations is done through the reformulation
\fref{characIm(R_eta)} of the equation satisfied by $\X^\eta$.
Moreover, we shall expand $\X^\eta$ in powers of $\eta^2$ up to the
first order.


\begin{definition}
Let $\overline{Q}$ be the operator on $ L^2_\MM$ defined by
\fref{qbar} and let $\qbarbar$ the one defined on
$\overline{L}^2_\MM$ by \fref{qbarbar}. Then $\qbar$ restricted to
$\overline{L}^2_\MM$ coincides with $\qbarbar$.
\end{definition}

We now list the properties of $\qbarbar$ which are completely inherited from those of $Q$.


\begin{proposition}[\emph{\textbf{Properties of $\overline{\overline{Q}}$}} ]\label{Qbar_prop1}

The operator $\qbarbar$ is a bounded operator on
$\overline{L}^2_\MM$ equipped with the $L^2_\MM$ scalar product. It
satisfies the following properties whose proof are immediate and are
left to the reader.
\begin{enumerate}
\item For all $f\in \overline{L}^2_\MM$, we have $\qbarbar(f) =\qbar(f) = \overline{Q(f)}$.
\item $-\qbarbar$ is a bounded, symmetric, non-negative operator on $\bar
L^2_\MM$.
\item The null set of $\qbarbar$  is given by
\begin{equation*}
N(\qbarbar)=\{n \MM(\v), n\in \RR\}
\end{equation*}
\item Let $\mathcal{P}:\overline{L}^2_\MM\rightarrow N(\qbarbar)$ be the
orthogonal projection on
 $N(\qbarbar)$, the following coercivity inequality holds,
\begin{equation*}
-\langle\qbarbar(f),f\rangle_\MM\geq \alpha_1\|f-\mathcal{P}(f)\|_\MM^2.
\end{equation*}
i.e $-\qbarbar$ is coercive on $N(\qbarbar)^\perp$.

\item The range of $\qbarbar$ in $\overline{L}^2_\MM$ is given by
\begin{equation*}
Im(\qbarbar)=N(\overline{Q})^\perp=\left\{\ds f\in
\overline{L}^2_\MM/\, \int_\v fd\v=0\right\}= Im(Q)\cap \overline{L}^2_\MM
\end{equation*}
\end{enumerate}

\end{proposition}


In the sequel, we shall also need the properties of $\qbar$
considered on the whole space $L^2_\MM$ (whereas $\qbarbar$ is
restricted to cylindrically symmetric functions). Some of the
properties of $\qbar$ are inherited from those of $Q$ but some of
these properties like the symmetry are lost. The following
proposition summarizes these properties.
\begin{proposition}[\emph{\textbf{Properties of $\overline{Q}$ }}]\label{Qbar_prop2}
The operator $\qbar$ is a bounded operator. It satisfies the following properties
\begin{enumerate}
\item The null set of $\overline{Q}$ is
\begin{equation*}
N(\overline{Q})= N( Q)=\{n \MM(v), n\in \RR\}.
\end{equation*}
\item The range of $\overline{Q}$ in $L^2_\MM$ is the following set
\begin{equation}
\label{solqbar}
Im(\overline{Q})=\left\{g\in L^2_\MM/ \int_{\RR^3}\frac{\nu g}{\bar
\nu}d\v=0\right\}.
\end{equation}
\item For each $g\in Im(\overline{Q})$, there exists a unique function $f\in
L^2_\MM$ satisfying  $\ds\int_{\RR^3}f d\v=0$ such that $\overline{Q}(f)=g$. The solution can be written
 $\ds f=\overline{f}+ \frac{\overline{g}- g}{\overline{\nu}}$ where the
isotropic part of f, $\overline{f} = {\mathcal{A}} (f)$, verifies
\begin{equation}\label{f_barequation}
\qbarbar(\overline{f})=\overline{g} + \overline{Q}^+(\frac{g-\overline{g}}{\overline{\nu}}),
\end{equation}
and where we have denoted $\qbar = \overline{Q}^+ - \overline{\nu}.$
In addition, there exists a constant $C>0$ independent of $g$ (and
$f$) such that
\begin{equation}\label{Q_bar_bound}
\|f\|_\MM\leq C\|g\|_\MM.
\end{equation}
\end{enumerate}

\end{proposition}
\begin{proof}
Let us begin with item 1. It is clear that if $f$ is  Maxwellian, it
is in the kernel of $\qbar$. Let now $f\in N(\overline{Q})$. We have
$\overline{Q}^+(f)-\overline{\nu} f=0 $. Since the Maxwellian is
cylindrically symmetric, it is readily seen that $\overline{Q}^+(f)$
is always cylindrically symmetric. Therefore, $f = {1\over
\overline{\nu} } \overline{Q}^+(f)$, is cylindrically symmetric and
thus
 $$\qbarbar (f) = \qbar(f) = 0.$$
In view of Proposition \ref{Qbar_prop1}, $f$ is proportional to the
maxwellian.

\medskip In order to prove items 2. and 3., we consider $g\in Im(\overline{Q})$, and  $f\in L^2_\MM$ such that
$\overline{Q}(f)=g$. Decomposing $f$ and $g$ as $f=\overline{f}+h_f$ (resp.
$g=\overline{g}+h_g$) where $\overline{h_f}=\overline{h_g}=0$, we have
$$\overline{Q}(\overline{f}) + \overline{Q}^+ (h_f)  - \overline{g} = h_g +\overline{\nu} h_f.$$
The left hand side of the above equality is cylindrically
symmetric, whereas cylindrical averages of the right hand side
vanish. Therefore, both are equal to zero. We deduce that
$$h_f = -\frac{h_g}{\overline{\nu}}=\frac{\overline{g}- g}{\overline{\nu}},$$
whereas $\overline{f}$ satisfies
$$\qbar(\overline{f}) = \overline{g}  - \overline{Q}^+ (h_f)$$
which can be rewritten
$\qbarbar(\overline{f}) = \overline{g} - \overline{Q}^+( \frac{\overline{g}- g}{\overline{\nu}})$
and also
$$\qbarbar(\overline{f} + \frac{\overline{g}}{\overline{\nu}}) =  \overline{Q}^+( \frac{ g}{\overline{\nu}}).$$
This equation has a solution if and only if the right hand side has
zero velocity average (see Proposition \ref{Qbar_prop1}). A
straightforward computation shows that this is equivalent to the
condition $\ds \int_{\RR^3}\frac{\nu g}{\bar \nu}d\v=0$. Let us
finally show estimate \eqref{Q_bar_bound}. Remark first that for any
$f\in L^2_\MM$ with $f=\overline f+ h_f$ and $\overline{h_f}=0$, we have
$\|f\|_\MM^2=\|\overline{f}\|_\MM^2+\|h_f\|_\MM^2$ and then
$\|\overline{f}\|_\MM \leq \|f\|_\MM$. Besides, equation (\ref{f_barequation}) is posed on
$\overline{L}^2_\MM$ then, from the coercivity of $-\qbarbar$ on
$N(\qbarbar)^\perp$ (see Proposition \ref{Qbar_prop1}), there exists
a constant $C_1>0$ such that
$$
\|\overline{f}\|_\MM\leq C_1(\|\overline{g}+\overline{Q^+}(\frac{g-\bar
g}{\overline{\nu}})\|_\MM) \leq C_1(\|g\|_\MM+\|Q^+(\frac{g-\bar
g}{\overline{\nu}})\|_\MM) \leq C_1' \|g\|_\MM.$$
In addition, since $\ds h_f=\frac{\overline{g}- g}{\overline{\nu}}$ and
with \eqref{crosssection}, there exists $C_2 >0$ such that $\|h_f\|_\MM\leq C_2 \|g\|_\MM$ and
estimate \eqref{Q_bar_bound} holds.

\end{proof}


\begin{definition}
\label{rotation} We define the gyration operator $\G$ by $\G(f) =
(\v\times \e_z) \cdot \nabla_\v f$, with domain $D$.
\end{definition}
The kernel of this operator is nothing but the set of cylindrically
symmetric functions. It is clear that its range is contained in the
set of $L^2_\MM$ functions which have zero cylindrical averages. The
following lemma shows that both sets coincide.

\begin{lemma}\label{gg}
Let  $g\in L^2_\MM$ with zero cylindrical average (${\mathcal{A}} (g) = 0$). Then, the set of solutions of the equation
\begin{equation}
\label{g}
\G (f) = g,
\end{equation}
in nonempty. Any solution $f$ can be written
$$f = f_1 +  {\mathcal{A}}_1 (g)$$
where $f_1$ si an arbitrary cylindrically symmetric function (which means that ${\mathcal{A}}(f_1) = f_1$) and the
average operator  ${\mathcal{A}}_1$ is defined by
\begin{equation}
\label{gammafunctions}
{\mathcal{A}}_1(g)={1\over \overline{\nu} (\v)}\int_0^{2\pi} g(\RC(\tau)v)\overline{\nu_\tau }d\tau =
{1\over 2\pi \overline{\nu}} \int_0^{2\pi} g(\RC(\tau)v)
\int_0^\tau \nu(\RC (s)v)\, ds\, d\tau.
\end{equation}
Moreover, we have
$${\mathcal{A}} (\nu {\mathcal{A}}_1 (g)) = 0.$$
\end{lemma}
\begin{proof}
Let us consider the problem
$$\G (f_\eta) + \eta^2 \nu f_\eta = g.$$
This problem has a unique solution $f_\eta$ defined by
$$\begin{array}{lll}
\ds f_\eta = {1\over \eta^2} L_\eta (g) &= &\ds\frac{1}{
C_\eta(v)}\int_0^{2\pi}g(\RC(\tau)v)
{e^{\eta^2\int_0^{\tau}\nu(\RC(s)v)ds}}d\tau \\[3mm]
&= &\ds \frac{1}{
C_\eta(v)}\int_0^{2\pi}g(\RC(\tau)v)
[{e^{\eta^2\int_0^{\tau}\nu(\RC(s)v)ds}}- 1]d\tau
 \end{array}
 $$
where
$$
\ds \ds C_\eta(v)=\ds e^{ \eta^2 \int_0^{2\pi}\nu(\RC(s)v)ds}-1 = e^{2\pi\eta^2 \overline{\nu}} - 1.
$$
Passing to the limit $\eta\to 0$, we find that $f_\eta$ converges
almost everywhere towards ${\mathcal{A}}_1(g).$ Moreover ${\mathcal{A}}(\eta^2
\nu f_\eta) = {\mathcal{A}}(g) = 0$, which proves after passing to the
limit $\eta\to 0$ that ${\mathcal{A}}(\nu{\mathcal{A}}_1(g)) = 0$. Of course,
any other solution of the equation $\G (f) = g$ is obtained by
adding an arbitrary cylindrically symmetric function.

\end{proof}

Let us now consider the problem
\begin{equation}
\label{ga}
\left\{
\begin{array}{lll}
\G(f) &= &g \\
{\mathcal{A}} (Q(f))& = &h.
\end{array}
\right.
\end{equation}
It is clear that the above system does not have a solution unless
${\mathcal{A}}(g) = 0$ and $h\in \overline{L}^2_\MM$ with $\int_{\RR^3}
h\, d\v = 0$. The proposition below shows that these conditions are
sufficient and gives the general solution of the problem
\begin{proposition}\label{SQ}
 Let $g$ and $h$ be in $L^2_\MM$ such that ${\mathcal{A}}(g) =
0$, ${\mathcal{A}} (h) = h$ and $\int h\, d\v = 0$. Then the problem
\fref{ga} has a unique solution up to a Maxwellian (the Maxwellian
is the only solution of the homogeneous problem). The unique
solution $f$ which has zero average is given by
$$f = f_1 + {\mathcal{A}}_1(g)$$
$$\qbar (f) =- {\overline{\nu}}{\mathcal{A}}_1(g) +h$$
with $f_1$ is an arbitrary cylindrically symmetric function.
Moreover, the mapping $(g,h) \mapsto f$ is linear and continuous
both on $L^2_\MM$ and $H^2_\MM$, i.e.
$$\|f\|_{L^2_\MM} \leq C (\|h\|_{L^2_\MM} + \|g\|_{L^2_\MM}),\quad \|f\|_{H^2_\MM} \leq C (\|h\|_{H^2_\MM} + \|g\|_{H^2_\MM}).$$
\end{proposition}

\begin{proof}
First of all, Lemma \ref{gg} shows that $f = f_1 + {\mathcal{A}}_1 (g)$.
Besides, we have
$$\overline{Q}(f) = {\mathcal{A}}(Q(f)) + {\mathcal{A} } (\nu f) - \overline{\nu} f.$$
Since $f_1$ is cylindrically symmetric, then ${\mathcal{A} } (\nu f_1) -
\overline{\nu} f_1 = 0$, which leads, thanks to the identity ${\mathcal{A}}(\nu {\mathcal{A}}_1(g)) = 0$, to $\qbar (f) =- {\overline{\nu}}{\mathcal{A}}_1(g) +h.$ It is now enough, in view of the solvability condition
\fref{solqbar} to check that $\int{\nu \over \overline{\nu}} (-
{\overline{\nu}}{\mathcal{A}}_1(g) +h) d\v = 0,$ which is readily seen.

\end{proof}

\subsection{Expansion of $\mathrm{X}^\eta_z$}
The idea is to make a Hilbert expansion of $\mathrm{X}^\eta_z$ in
powers of $\eta^2$ :
$$\mathrm{X}^\eta_z = \sum_{i} \eta^{2i} \mathrm{X}^{(i)}_z.
$$
We have the following equations

$$
\begin{array}{lll}
\G(\mathrm{X}^{(0)}_z) &= &0,\\[3mm]
\G(\mathrm{X}^{(1)}_z) &= &Q (\mathrm{X}^{(0)}_z) +v_z \MM,\\[3mm]
\G(\mathrm{X}^{(i+1)}_z) &= &Q (\mathrm{X}^{(i)}_z), \quad i\geq 1,
\end{array}
$$
where we have taken advantage of the fact that ${\mathcal{A}} (v_z \MM) = v_z \MM.$
The solvability conditions become
$${\mathcal{A}}(Q (\mathrm{X}^{(0)}_z) + v_z \MM) = 0$$
$${\mathcal{A}}(Q (\mathrm{X}^{(i)}_z)) = 0, \quad i\geq 1.$$
Therefore, the terms of the expansion can be successively computed
by solving the problems
$$
\begin{array}{l}
\left\{
\begin{array}{lll}
\G(\mathrm{X}^{(0)}_z) &= &0 \\
{\mathcal{A}} (Q(\mathrm{X}^{(0)}_z))& = &-v_z \MM.
\end{array}
\right.
\\[8mm]
\left\{
\begin{array}{lll}
\G(\mathrm{X}^{(1)}_z) &= &Q (\mathrm{X}^{(0)}_z) + v_z \MM. \\
{\mathcal{A}} (Q(\mathrm{X}^{(1)}_z))& = &0.
\end{array}
\right.
\\[8mm]
\left\{
\begin{array}{lll}
\G(\mathrm{X}^{(i+1)}_z) &= &Q (\mathrm{X}^{(i)}_z) \\
{\mathcal{A}} (Q(\mathrm{X}^{(i+1)}_z))& = &0.
\end{array}
\right.
\end{array}
$$
 From  Proposition \ref{SQ}, we deduce that $\mathrm{X}_z^{(i)}$ are
uniquely determined by the above equations and are in $H^2_\MM$.
\begin{proposition}\label{proposition_chieta||}
Let $\X^\eta$ be defined by \fref{chi_eta}. The function
$\mathrm{X}^{(0)}_z$ is uniquely determined by
\begin{equation}\label{chi_z^0}
-\qbarbar(\mathrm{X}_z^{(0)})=v_z\MM,\quad \int_{\RR^3}\mathrm{X}_z^{(0)}(\v)d\v=0
\end{equation}
and $\mathrm{X}_z^{(1)}$ is the unique solution in $L^2_\MM$ of
\begin{equation}\label{chi_z^1}
-\overline{Q}(\mathrm{X}_z^{(1)})=\bar \nu{\mathcal{A}}_1(Q(\mathrm{X}_z^{(0)}) + v_z\MM)= \bar \nu{\mathcal{A}}_1(Q(\mathrm{X}_z^{(0)}) -\qbarbar(\mathrm{X}_z^{(0)})),\quad
\int_{\RR^3}\mathrm{X}_z^{(1)}(\v)d\v=0
\end{equation}
\noindent where ${\mathcal{A}}_1$ is the averaging operator defined by
\eqref{gammafunctions}.
Then the following expansion holds in the strong topology of $\HC^2_\MM$
\begin{equation}\label{chi_eta^z_expansion}
\mathrm{X}^\eta_z = \mathrm{X}_z^{(0)}+ \eta ^2 \mathrm{X}_z^{(1)} + {\mathcal{O}} (\eta^4).
\end{equation}
\end{proposition}

\begin{proof}
 We write
$$\mathrm{X}^\eta_z = \mathrm{X}_z^{(0)}+ \eta^2 \mathrm{X}_z^{(1)} + \eta^4  r^\eta.$$
From the definitions of $\mathrm{X}_z^{(i)}$, it is readily seen
that
$$\G (r^\eta) =\eta^2 Q(r^\eta) + Q( \mathrm{X}_z^{(1)}) $$
$${\mathcal{A}} (Q(r^\eta )) = 0.$$
Hence, we have that
$$\overline{Q}(r^\eta) = - \overline{\nu} {\mathcal{A}}_1 (\eta^2Q(r^\eta) +  Q(\mathrm{X}_z^{(1)})).$$
Using estimate \eqref{Q_bar_bound}, we obtain for some constant $C
>0$ independent of $\eta$
$$\|r^\eta\|_{L^2_\MM} \leq C (\eta^2 \|r^\eta\|_{L^2_\MM} +  \|\mathrm{X}_z^{(1)}\|_{L^2_\MM})$$
which shows the boundedness of $r^\eta$ in $L^2_\MM$.
Differentiating the above identity with respect to the velocity
variable, on can show simply the boundedness of $ r^\eta$ in
$\HC^2_\MM$.
\end{proof}


\subsection{Expansion of $\X_\eta^\perp$}
In this subsection, we expand the orthogonal part of $\X^\eta$. We
follow the same steps as for $\mathrm{X}^\eta_z$ and write the
equations satisfied by the formal expansion of $\X^\eta_\perp$.
Namely, we have

$$
\begin{array}{lll}
\G(\X^{(0)}_\perp) &= &\v_\perp \MM, \\[3mm]
\G(\X^{(1)}_\perp) &= &Q (\X^{(0)}_\perp),\\[3mm]
$$\G(\X^{(i+1)}_\perp) &=& Q (\X^{(i)}_\perp), \quad i\geq 1.
\end{array}
$$
As for $\mathrm{X}^\eta_z$, the above system can be reformulated as follows
$$
\begin{array}{l}
\left\{
\begin{array}{lll}
\G(\X^{(0)}_\perp) &= &\v_\perp \MM \\
{\mathcal{A}} (Q(\X^{(0)}_\perp))& = &0.
\end{array}
\right.
\\[8mm]
\left\{
\begin{array}{lll}
\G(\X^{(1)}_\perp) &= &Q (\X^{(0)}_\perp) \\
{\mathcal{A}} (Q(\X^{(1)}_\perp))& = &0.
\end{array}
\right.
\\[8mm]
\left\{
\begin{array}{lll}
\G(\X^{(i+1)}_\perp) &= &Q (\X^{(i)}_\perp) \\
{\mathcal{A}} (Q(\X^{(i+1)}_\perp))& = &0.
\end{array}
\right.
\end{array}
$$

\begin{proposition}\label{proposition_chietaperp}
The expansion
$$\X^\eta_\perp  = \X^{(0)}_\perp + \eta^2 \X^{(1)}_\perp + {\mathcal{O}} (\eta^4),$$
holds true in $\HC^2_\MM$.
Morever, $\X^{(0)}_\perp$ satisfies
$$\X^{(0)}_\perp - {\mathcal{A}} (\X^{(0)}_\perp)= -{\mathcal{I}} (\v_\perp \MM),$$
with $\II=\left(
                                    \begin{array}{cc}
                                      0 &1 \\
                                      -1& 0 \\
                                    \end{array}
                                  \right)
.$
\end{proposition}
The proof of this proposition is identical to that of Proposition
\ref{proposition_chieta||} and is skipped. The only thing which has
to be checked is the expression of the anisotropic part of
$\X^{(0)}_\perp$ which can be deduced from the identities
$$-\v_\perp \MM =   \G ({\mathcal{I}} (\v_\perp \MM)), \quad {\mathcal{A}}({\mathcal{I}} (\v_\perp \MM)) = 0.$$


\section{Proof of the main theorems}

We begin this section by making precise the definition of weak
solutions of the Boltzmann equation and give quite standard a priori
estimates on this solution. For the moment, $\ep$ and $\eta$ are
arbitrary and all the constants are, unless specified, independent
of these parameters.

\begin{definition}
 Let $T \in \RR^*_+$, a function $f^{\ep\eta}\in L^1_{loc}([0,T]\times \RR^6)$ is called weak solution of
\eqref{boltzmann_eps.eta} if it satisfies

\begin{multline}\label{weaksolution}
-\int_0^T\int_{\RR^6} f^{\ep\eta}\left[ \pa_t \psi  +  {\T_z\psi \over \ep} + { \T_\perp \psi\ \over \ep\eta }
\right] dt d\r d\v-\frac{1}{\ep^2\eta^2}\int_0^T\int_{\RR^6}
f^{\ep\eta} \G(\psi) dtd\r d\v\\=
\frac{1}{\ep^2}\int_0^T\int_{\RR^6} Q(f^{\ep\eta})\psi dt d\r d\v +
\int_{\RR^6}f_0(x,v)\psi(0,x,v)d\r d\v.
\end{multline}

\noindent for all $\psi \in C^1_c([0,T)\times \RR^6)$.
\end{definition}

\begin{theorem}
\label{existence}
 Assume that $f_0\in L^2_{\MM_V}(d\r d\v)$. Then, $\forall \ep,\eta>0$ and
$\forall T \in \RR^*_+$, there exists a unique weak solution
$f^{\ep\eta}\in C^0([0,T], L^2_{\MM_V}(d\r d\v))$ of
\eqref{boltzmann_eps.eta}. In addition, the following mass
conservation equation holds
\begin{equation}\label{continuity.eps.eta}
\partial_t\rho^{\ep\eta}+\mathrm{div}_\r \mathbf{J}^{\ep\eta}=0
\end{equation}
where
\begin{equation}\label{density}
\rho^{\ep\eta}(t,x)=\int_{\RR^3}f^{\ep\eta}d\v, \quad
\mathrm{J}^{\ep\eta}_z =\frac{1}{\ep}\int_{\RR^3}v_z
f^{\ep\eta}d\v, \quad
\mathbf{J}^{\ep\eta}_\perp=\frac{1}{\ep\eta}\int_{\RR^3}\v_\perp
f^{\ep\eta}d\v,
\end{equation}
and there
exists a constant $C>0$ independent of $\ep$ and $\eta$ such that
\begin{equation}\label{fbound}
 \|f^{\ep\eta}(t)\|_{L^2_\MM(d\r d\v)}\leq C \qquad \forall t\in[0,T],
\end{equation}
\begin{equation}\label{ecarttype}
\|f^{\ep\eta}-\rho^{\ep\eta}\MM\|_{L^2([0,T],L^2_\MM(d\r d\v))}\leq C \ep
\end{equation}
\end{theorem}
An immediate consequence of the integrability properties of the
function $f^{\ep\eta}$, with respect to $\v$ is that the weak
formulation \fref{weaksolution} is satisfied by test functions which
are not necessarily compactly supported in velocity. More precisely,
we have
\begin{corollary}
Let $f^{\ep\eta}$ be the solution exhibited in Theorem
\ref{existence}. Then the weak formulation \fref{weaksolution} is
satisfied by test functions lying in the set
\begin{multline}\label{testfunctions} \mathcal{T}=\{\psi(t,\r,\v)\in
C^1([0,T)\times \RR^{6}) \mbox{ compactly supported w.r.t. $t$ and
$\r$} \mbox{ and } \\ \qquad \exists n\in \NN, C\in\RR_+, \quad
|\nabla_{t,\r,\v} \psi(t,\r,\v)|\leq  C (1+|\v|)^n \}.
\end{multline}
\end{corollary}
We shall not prove this corollary which is straightforward. The
proof of Theorem \ref{existence} is also classical (see for example
\cite{poupaud3}). We show here how the a priori estimates
\fref{fbound} and \fref{ecarttype} can be obtained. We simply
multiply the Boltzmann equation \fref{boltzmann_eps.eta} by
${f^{\ep\eta}\over \MM_V}$ and integrate with respect to $\r$ and
$\v$. Straightforward computations show that this leads to the
entropy inequality

\begin{align*}&{d\over dt} \int {(f^{\ep\eta})^2 \over
2\MM_V}\, d\r d\v -\int \pa_t V {(f^{\ep\eta})^2 \over 2\MM_V}\, d\r
d\v  \\
=& {1\over \ep^{2}} \int Q^\eta(f^{\ep\eta}) {f^{\ep\eta}\over
\MM_V}\, d\r d\v \leq - {\alpha_1\over \ep^2} \|f^{\ep\eta} -
\rho^{\ep\eta} \MM\|_{L^2_{\MM_V}}^2.\end{align*} Thanks to the
boundedness of the potential $V$ and its derivatives, estimate
\fref{fbound} follows by applying Gronwall Lemma and then the bound
\fref{ecarttype} follows immediately.


\subsection{Proof of Theorem \ref{main_th.2}}
In this section, we are interested in the limit $\ep \to 0$, while
$\eta$ is fixed. To this aim, we follow the usual approach based on
the moment method. Namely, thanks to estimates \eqref{fbound} and
\eqref{ecarttype}, the sequences $f^{\ep\eta}, (\rho^{\ep\eta}),
(J^{\ep\eta}_z)$ and $(\mathbf{J}^{\ep\eta}_\perp)$ are bounded
with respect to $\ep$. One can find $\rho_\eta(t,x)\in
L^2((0,T)\times \RR^3)$ such that $\rho^{\ep\eta}$ converges weakly
to $\rho_\eta$ in $L^2((0,T)\times \RR^3)$ and
$f^{\ep\eta}\rightharpoonup \rho_\eta \MM$ in
$L^\infty(0,T;L^2_\MM(\RR^6))$. Similarly, $\mathbf{J}^{\ep\eta}$
converges weakly to $\mathbf{J}^\eta$ in $L^2((0,T)\times \RR^3)$.
We have of course, the mass conservation equation
$$\pa_t\rho_{\eta} +\mathrm{div}_\r \mathbf{J}^\eta= 0.$$
 Now the only thing left to show is to identify $\mathbf{J}^\eta$.
The way to proceed for identifying $\mathbf{J}^\eta$ is standard
(see for instance \cite{bar,bensouss,sentis, dautray, poupaud3,
Deg_Gou_Pou, Babd_Deg,Babd_Deg_Gen,Babd_Desv_Gen}). The Boltzmann
equation can be rewritten
\begin{equation}\label{Boltz2}
\ds\ep \frac{\partial
 f^{\varepsilon\eta}}{\partial t}+ {\T_z f^{\varepsilon\eta}}
 +\frac{\T_\perp
f^{\ep\eta}}{\eta} =\frac{Q^{\eta}(f_{\varepsilon\eta})}{\varepsilon}
\end{equation}
Let now $\X^{\eta*}$ be the solution of
$$Q^{\eta*} (\X^{\eta*} )= - \left(\begin{array}{c} \ds {1\over \eta^2}\v_\perp\\ v_z\end{array}\right)\MM(\v)\quad \mbox{and }\quad
\int_{\RR^3}\X^{\eta*}(\v)d\v=0,
$$
where $Q^{\eta*}:=(Q^{\eta})^*$ is the adjoint in $L^2_\MM$ of
$Q^\eta$. It is readily seen that $Q^{\eta*}$ has the same
expression as $Q^\eta$ except that $\eta^2$ is replaced by
$-\eta^2$. We then have the following result.
\begin{corollary}
\label{chi*} The following expansion holds in $\HC^2_\MM$ $$\X^{\eta
*}_\perp = - \X^{(0)}_\perp + \eta^2 \X^{(1)}_\perp + {\mathcal{O}}(\eta^4)\quad; \quad X^{\eta*}_z = X^{(0)}_z - \eta^2 X^{(1)}_z +
{\mathcal{O}}(\eta^4).
$$
\end{corollary}
We multiply \eqref{Boltz2} by $\ds {X_z^{\eta*} \over \MM}$ and
integrate w.r.t. $\v$. The right hand side of the obtained equation
is nothing but
$$\frac{1}{\ep}\int{Q^{\eta} (f^{\ep\eta}) X_z^{\eta*} \over \MM} d\v =\frac{1}{\ep} \int
{Q^{\eta *} (X_z^{\eta*} ) f^{\ep\eta} \over \MM} d\v = -
J^{\ep\eta}_z. $$ Therefore $$- J^{\ep\eta}_z = \ep \pa_t \int
{f^{\ep\eta} X_z^{\eta*} \over \MM}\, d\v + \int {X_z^{\eta*} \over
\MM} \T_z f^{\ep\eta} d\v + {1\over \eta} \int {X_z^{\eta*} \over
\MM}  \T_\perp f^{\ep\eta} d\v.
$$
Thanks to the bounds \fref{fbound} and
\fref{polyn.increas.chi_eta_1}, we can pass to the limit in the
above equation as $\ep$ tends to zero and obtain $$- J^{\eta}_z =
\int {X_z^{\eta*} \over \MM} \T_z (\rho_{\eta} \MM) d\v + {1\over
\eta} \int {X_z^{\eta*} \over \MM}  \T_\perp (\rho_{\eta}\MM) d\v,
$$
which leads to the expression $$-J^{\eta}_z = [\int X_z^{\eta*}
v_zd\v] (\pa_z \rho_\eta + \rho_\eta \pa_z V) + {1\over \eta} [\int
X_z^{\eta*} \v_\perp d\v]\cdot (\nabla_{\r_\perp} \rho_\eta+
\rho_\eta \nabla_{\r_\perp} V).$$ Multiplying the Boltzmann equation
\eqref{Boltz2} by $\ds \frac{\eta \X^{\eta*}_\perp}{\MM}$ and
proceeding analogously, we find out that
$$-\mathbf{J}^\eta_\perp = \eta [\int \X_\perp^{\eta*} v_zd\v] (\pa_z \rho_\eta + \rho_\eta \pa_z V)
+ [\int \X_\perp^{\eta*}\otimes\v_\perp  d\v] (\nabla_{\r_\perp}
\rho_\eta + \rho_\eta \nabla_{\r_\perp} V).$$ Using the property
$\int X_z^{\eta*} v_zd\v = \int v_z X_z^\eta d\v$ as well as $\int
X_z^{\eta*} \v_\perp d\v = \eta^2 \int v_z \X^{\eta}_\perp d\v$ and
analogous formulae for the terms appearing in the expression of
$\mathbf{J}^\eta_\perp$, we finally obtain the expression
$$\mathbf{J}^\eta = - \DD^\eta (\nabla_\r \rho_\eta + \rho_\eta \nabla_r V)$$
where
$$\DD^\eta = \int \left(\begin{array}{c} {\eta} \X^{\eta*}_\perp \\ X^{\eta*}_z
\end{array}\right)\otimes \left( \begin{array}{c} {1\over \eta} \v_\perp \\ v_z
\end{array}\right)  d\v = \int \left(\begin{array}{c} {1\over \eta} \v_\perp \\ v_z
\end{array}\right)\otimes\left( \begin{array}{c} \eta \X^\eta_\perp \\ X^\eta_z
\end{array}\right)   d\v.
$$
One can then expand $\DD^\eta$ as

$$\DD^\eta = \DD^\eta_\perp + \DD^\eta_z + \eta \DD^\eta_{z\perp}$$
where
$$\DD^\eta_\perp = \int \left(\begin{array}{c}  \X^{\eta*}_\perp \\ 0
\end{array}\right)\otimes\left( \begin{array}{c} \v_\perp \\ 0
\end{array}\right)   d\v = \int\left(\begin{array}{c}  \v_\perp  \\ 0
\end{array}\right) \otimes\left( \begin{array}{c} \X^{\eta}_\perp \\ 0
\end{array}\right)   d\v,
$$
$$D^\eta_z =
\int \left(\begin{array}{c} \mathbf{0} \\ X^{\eta*}_z
\end{array}\right)\otimes \left(\begin{array}{c} \mathbf{0} \\ v_z
\end{array}\right)   d\v =
\int \left(\begin{array}{c} \mathbf{0} \\ v_z
\end{array}\right)\otimes\left( \begin{array}{c} \mathbf{0} \\ X^{\eta}_z
\end{array}\right)   d\v,
$$
$$
\DD^\eta_{z\perp} = \int \left(\begin{array}{c} \mathbf{0} \\ v_z
\end{array}\right)\otimes\left( \begin{array}{c} \X^\eta_\perp \\ 0
\end{array}\right)   d\v
+ \int\left(\begin{array}{c}  \X^{\eta*}_\perp \\ 0
\end{array}\right)\otimes \left( \begin{array}{c} \mathbf{0} \\ v_z
\end{array}\right)   d\v.
$$

\subsection{Proof of Proposition \ref{chieta}}
Let us expand first $\DD^\eta_\perp$. Following the expansion of
$\X^{\eta*}_\perp$ and that of  $\X^{\eta}_\perp$, we obtain that
$$\DD^\eta_\perp = \DD^{(0)}_\perp + \eta^2 \DD^{(1)}_\perp +{\mathcal{O}}(\eta^4),$$
where

$$\DD^{(0)}_\perp =
- \int \left(\begin{array}{c}  \X^{(0)}_\perp \\ 0
\end{array}\right)\otimes \left( \begin{array}{c} \v_\perp \\ 0
\end{array}\right)  d\v
=
 \int \left(\begin{array}{c}  \v_\perp \\ 0
\end{array}\right)\otimes\left( \begin{array}{c} \X^{(0)}_\perp \\ 0
\end{array}\right)   d\v
$$
and
$$\DD^{(1)}_\perp =
\int \left(\begin{array}{c}  \X^{(1)}_\perp \\ 0
\end{array}\right)\otimes \left( \begin{array}{c} \v_\perp \\ 0
\end{array}\right)  d\v
=
 \int \left(\begin{array}{c}  \v_\perp \\ 0
\end{array}\right)\otimes \left( \begin{array}{c} \X^{(1)}_\perp \\ 0
\end{array}\right)  d\v.
$$
It is clear that $\DD^{(0)}_\perp$ is antisymmetric while
$\DD^{(1)}_\perp$ is symmetric and we have
$$\DD^{(0)}_\perp =  - \int \left(\begin{array}{c}  \X^{(0)}_\perp  - {\mathcal{A}}(\X^{(0)}_\perp)\\ 0
\end{array}\right)\otimes\left( \begin{array}{c} \v_\perp \\ 0
\end{array}\right)   d\v.
$$
Since $\X^{(0)}_\perp  - {\mathcal{A}}(\X^{(0)}_\perp) = -{\mathcal{I}}
(\v_\perp \MM)$, we have

\begin{align*}\DD_\perp^{(0)} =&\left(\begin{array}{cc}\ds \int_{\RR^3}{\mathcal{I}}(\v_\perp\MM)\otimes\v_\perp dv & \textbf{0}\\
\textbf{0} & 0 \end{array}\right)=\left(\begin{array}{cc}\ds{\mathcal{I}}\cdot\int_{\RR^3}(\v_\perp\otimes \v_\perp)\MM dv & \textbf{0} \\
\textbf{0}& 0 \end{array}\right)\\
=&\left(\begin{array}{ccc} 0 & 1& 0\\
-1 &0 &0\\ 0& 0& 0 \end{array}\right)\end{align*} where we have
noticed that $\ds \int_{\RR^3} \v_\perp \otimes \v_\perp \MM d\v$ is
nothing but the identity matrix  on $\RR^2$. Let us now expand
$\DD^\eta_z$. We have
$$\DD^{\eta}_z = \DD^{(0)}_z + \eta^2 \DD^{(1)}_z + {\mathcal{O}}(\eta^4)$$ where
$$\DD^{(0)}_z = D_z \left(\begin{array}{ccc}
0&0&0\\[3mm]
0&0&0\\[3mm]
0&0&1
\end{array}\right)
$$
which $D_z = \int v_z \X^{(0)}_z\, d\v$. Comparing the expansions of
$X^\eta_z$and $X^{\eta*}_z$, we deduce that $\DD^{(1)}_z$ is equal
to zero. Let us now expand $\DD^\eta_{z\perp}$. Of course we have
$$\DD^\eta_{z\perp} = \DD^{(0)}_{z\perp} + \eta^2 \DD^{(1)}_{z\perp} +{\mathcal{O}}(\eta^4),$$
where
$$
\DD^{(0)}_{z\perp} = \int \left(\begin{array}{c} \mathbf{0} \\ v_z
\end{array}\right)\otimes\left( \begin{array}{c} \X^{(0)}_\perp \\ 0
\end{array}\right)   d\v
- \int \left(\begin{array}{c}  \X^{(0)}_\perp \\ 0
\end{array}\right)\otimes\left( \begin{array}{c} \mathbf{0} \\ v_z
\end{array}\right)   d\v
$$
and
$$
\DD^{(1)}_{z\perp} = \int \left(\begin{array}{c} \mathbf{0} \\ v_z
\end{array}\right)\otimes\left( \begin{array}{c} \X^{(1)}_\perp \\ 0
\end{array}\right)   d\v
+ \int\left(\begin{array}{c}  \X^{(1)}_\perp \\ 0
\end{array}\right)\otimes \left( \begin{array}{c} \mathbf{0} \\ v_z
\end{array}\right)   d\v.
$$
We can now make explicit the expansion of $\DD^\eta = \DD^{\eta}_s +
\DD^\eta_{as}$ where the indices $s$ and $as$ stand for the
symmetric and antisymmetric parts
$$\DD^{\eta}_s = \left(\begin{array}{cc}\eta^2 \int \v_\perp \otimes \X^{(1)}_\perp d\v & 0\\[3mm]
0 & D_z\end{array} \right) + \eta^3 \DD_{z\perp}^{(1)} +{\mathcal{O}}(\eta^4)$$
$$\DD^{\eta}_{as} = \left(\begin{array}{cc}{\mathcal{I}} & \bf{0}\\[3mm]
\bf{0} & 0\end{array}\right) + \eta \DD^{(0)}_{z\perp}+{\mathcal{O}}(\eta^4).
$$

\subsection{Proof of Theorem \ref{main_th.1}}

\begin{lemma} The density
$\rho^{\ep\eta}$  and the parallel current
$J^{\ep\eta}_z$ \eqref{density} are bounded in
$L^2_{t,x}([0,T]\times \RR^3)$ with respect to $\ep$ and $\eta$. In
addition, there exist $\rho$ and $\mathrm{J}_z$ in $L^2([0,T]\times
\RR^3)$ such that, up to extraction of subsequences, we have
\begin{eqnarray*}
&&\rho^{\ep\eta}\rightharpoonup \rho \quad \mbox{weakly in } L^2_{t,x}\\
&&f^{\ep\eta}\rightharpoonup \rho \MM \quad \mbox{in } L^\infty(
[0,T],L^2_\MM(d\r d\v))\\
&&\mathrm{J}_z^{\ep\eta}\rightharpoonup \mathrm{J}_z \quad \mbox{in }
L^2_{t,x}.
\end{eqnarray*}
\end{lemma}

\noindent However, estimates \eqref{ecarttype} and \eqref{fbound} do
not give a bound to $\mathbf{J}^{\ep\eta}_\perp$ and are not sufficient to
close the limit equation or to find relations between the
corresponding limits of the current and the density. To deal with
this lack of compactness, we have to filter out the
oscillations generated by the magnetic field in the orthogonal
direction.
\begin{proposition}
Let $\mathbf{J}^{\ep\eta}_\perp$ be the perpendicular part of the current
given by \eqref{density}. Then,
\begin{equation*}
\mathbf{J}^{\ep\eta}_\perp\rightharpoonup \mathbf{J}_\perp \quad \mbox{in }
D'([0,T)\times \RR^3)
\end{equation*}
with
\begin{equation*}
\mathbf{J}_\perp= (\rho \mathbf{E}_\perp  - \nabla_{\r_\perp} \rho ) \times \e_z
\end{equation*}
\end{proposition}
\begin{proof}

\noindent Multiplying \eqref{boltzmann_eps.eta} by $\ds \frac{\eta
\X^{\eta*}_\perp}{\MM}$ and integrating with respect to $\v$, we get
\begin{equation*}
-\mathbf{J}^{\ep\eta}_\perp={\ep}{\eta}\partial_t\int_{\RR^3}\frac{f^{\ep\eta}\X^{\eta*}_\perp}{\MM}d\v
+\int_{\RR^3}\T_\perp
f^{\ep\eta}\frac{\X^{\eta*}_\perp}{\MM}d\v+{\eta}\int_{\RR^3}\T_z
f^{\ep\eta}\frac{\X^{\eta*}_\perp}{\MM}d\v.
\end{equation*}

More precisely, by taking a test function $\ds \frac{\eta
\X^{\eta*}_\perp}{\MM}\phi(t,x)$ in the weak formulation
\eqref{weaksolution}, where  $\phi(t,x)\in C^1_c([0,T)\times \RR^3)$
we have

\begin{eqnarray*}
&&\nonumber\int_0^T\int_{\RR^3}
\mathbf{J}^{\ep\eta}_\perp(t,x)\phi(t,x)dtdx\\
&=&
{\ep}{\eta}\int_0^T\int_{\RR^6}f^{\ep\eta}\frac{\X^{\eta*}_\perp}{\MM}\partial_t
\phi dtd\r d\v+\ \int_0^T\int_{\RR^6}f^{\ep\eta}\T_\perp(
\frac{\X^{\eta*}_\perp\phi}{\MM}) dtd\r d\v\\ \nonumber &&+
{\eta}\int_0^T\int_{\RR^6}f^{\ep\eta}
\T_z(\frac{\X^{\eta*}_\perp\phi}{\MM}) dtd\r d\v + {\ep}{\eta}
\int_{\RR^6} f_0(x,v)\frac{\X^{\eta*}_\perp}{\MM}\phi(0,x) dx
d\v.\\
\end{eqnarray*}

\noindent In view of Corollary \ref{chi*}, $ \X^{\eta*}_\perp
\rightarrow -\X^{(0)}_\perp$ strongly in $\HC^2_\MM$ while
$(f^{\ep\eta})$ converges to $\rho \MM$ weakly  $L^2(0,T,
L^2_\MM(d\r d\v))$. Therefore, we can pass to the limit in the above
formula and get
\begin{eqnarray*}
\lim_{\ep,\eta\to 0} \int_0^T\int_{\RR^3}
\mathbf{J}^{\ep\eta}_\perp(t,x)\phi(t,x)dtdx
&=&-\int_0^T\int_{\RR^6}\rho\MM
\T_\perp(\frac{\X^{(0)}_\perp\phi}{\MM})d\r d\v dt\\
&=&\int_0^T\int_{\RR^3} \rho {\mathcal{I}}
(\nabla_{\r_\perp}\phi-\nabla_{\r_\perp} V \phi)d\r \, dt,
\end{eqnarray*}
where
\begin{equation*}
{\mathcal{I}} = \left(\begin{array}{cc}
0 & 1\\[3mm]
-1&0
\end{array}\right) = - \int_{\RR^3}\X^{(0)}_\perp\otimes\v_\perp  \, d\v.
\end{equation*} This concludes the proof of the proposition.

\end{proof}

The only thing left to do is to give the expression of the limiting
parallel current. To this aim, we use the test function $\ds
\frac{\mathrm{X}^{\eta*}_z}{\MM}\phi(t,x)$ in the weak formulation
where $\phi$ is a regular compactly supported function

\begin{eqnarray}\label{J_parallelweak}
&&\nonumber\int_0^T\int_{\RR^3} J^{\ep\eta}_z(t,x)\phi(t,x)dtdx\\
\nonumber&=& \int_0^T\int_{\RR^6}f^{\ep\eta}
\T_z(\frac{\mathrm{X}^{\eta*}_z\phi}{\MM}) dtd\r d\v + \ep
\int_{\RR^6} f_0(x,v)\frac{\mathrm{X}^{\eta*}_z}{\MM}\phi(0,x) dx
d\v.\\ \nonumber
&&+\ep\int_0^T\int_{\RR^6}f^{\ep\eta}\frac{\mathrm{X}^{\eta*}_z}{\MM}\partial_t
\phi dtd\r
d\v+\frac{1}{\eta}\int_0^T\int_{\RR^6}f^{\ep\eta}\T_\perp(
\frac{\mathrm{X}^{\eta*}_z\phi}{\MM}) dtd\r d\v.
\\
\end{eqnarray}
The first term of the right hand side converges as $\ep$ and $\eta$ tend to zero towards
$$\int_0^T\int_{\RR^6} \rho \MM
\T_z(\frac{\mathrm{X}^{(0)}_z\phi}{\MM}) dtd\r d\v =
\int_0^T\int_{\RR^3} \rho D_z(\pa_z \phi - \phi \pa_z V)\, dt d\r.$$
The second and third terms obviously tend to zero. The following
lemma shows that the last term of \fref{J_parallelweak} converges to
zero.

\begin{lemma}
Let $\phi (t,\r)$ be a compactly supported $C^2$ function. Then, we have
\begin{equation*}
\lim_{\ep,\eta\to 0}
\frac{1}{\eta}\int_0^T\int_{\RR^6}f^{\ep\eta}\T_\perp(
\frac{\mathrm{X}^{\eta*}_z\phi}{\MM}) dtd\r d\v = 0
\end{equation*}
\end{lemma}

\begin{proof}
Since  we have $\mathrm{X}^{\eta*}_z = \mathrm{X}^{(0)}_z +{\mathcal{O}}(\eta^2)$ in $\HC^2_\MM$ strong, we can replace
$\mathrm{X}^{\eta*}_z$ by $\mathrm{X}^{(0)}_z$. Now, we remark that
$\mathrm{X}^{(0)}_z$ is cylindrically symmetric so that
$${\mathcal{A}} (\T_\perp(
\frac{\mathrm{X}^{(0)}_z\phi}{\MM})) =0.$$
Therefore, one can find a test function $\psi(t,x,v)$ in the set ${\mathcal{T}}$ such that
$$\T_\perp(
\frac{\mathrm{X}^{(0)}_z\phi}{\MM}) = \G(\psi).$$
The fact that $\psi$ lies in ${\mathcal{T}}$ comes from the formula
$\psi = {\mathcal{A}}_1 (\T_\perp(
\frac{\mathrm{X}^{(0)}_z\phi}{\MM})).$
Using this function $\psi$ in the weak formulation, we obtain
$$\frac{1}{\eta}\int_0^T\int_{\RR^6}f^{\ep\eta}\T_\perp(
\frac{\mathrm{X}^{(0)}_z\phi}{\MM}) dtd\r d\v =
\frac{1}{\eta}\int_0^T\int_{\RR^6}f^{\ep\eta} \G(\psi) dtd\r d\v.$$
Using the weak formulation \fref{weaksolution}, the right hand side
of the above identity can be immediately estimated as ${\mathcal{O}}
(\ep)$, which tends to zero as $\ep$ and $\eta$ tend to zero.

 \end{proof}

 \section{Concluding remarks}
We have proven in this paper that the diffusion limit of the
Boltzmann equation with an ultrastrong magnetic field leads to a
diffusion equation in the parallel direction and a guiding center
motion in the orthogonal one. The proof has been done for the linear
Boltzmann equation and in the case of a constant magnetic field and
involved the analysis of the  operator involving collisions and
gyrations around the electrostatic field. If the magnetic field has
a constant direction but smoothly varies in position and time while
staying away from zero, the analysis can be carried out without
difficulty and the results can be generalized. A more difficult
problem appear if one couples the Boltzmann equation to the scaled
Poisson equation or when the collision operator has nonlinear
features. In this case, the method used in this paper might not be
sufficient and the use of double scale limits might be necessary as
it is for gyrokinetic limits.

\medskip

\section*{Acknowledgements} The authors acknowledge support
by the ACI Nouvelles Interfaces des Math\'ematiques No. ACINIM
176-2004 entitled ``MOQUA" and funded by the French ministry of
research, the ACI Jeunes chercheurs no. JC1035 ``Mod\`eles
dispersifs vectoriels pour le transport \`a l'\'echelle
nanom\'etrique" as well as the projet  No.  BLAN07-2 212988 entitled
``QUATRAIN" and funded by the Agence Nationale de la Recherche.


\medskip
Received April 2008; revised June 2008.

\medskip

\end{document}